\documentclass[a4paper,10pt]{article}

\usepackage[utf8]{inputenc}

\usepackage{aky-layout2}
\usepackage{aky-math}
\toggletrue{akynochapters}

\newcommand{\op}{\mathrm{op}}

\newcommand{\A}{\sA}

\renewcommand{\C}{\mathbf{C}}

\newcommand{\SmProj}{\mathbf{SmProj}}

\renewcommand{\AA}{{\sA^\bull}}

\newcommand{\Ddg}{{\underline{\D}}}

\renewcommand{\H}{\mathrm{H}}
\newcommand{\DGqeq}{\DGCat[\sW_{\mathrm{dg}}^{-1}]} 

\renewcommand{\L}{\bL}

\newcommand{\DerCorr}{\SmProj^\mathrm{cor}}

\newcommand{\ChMot}{\mathbf{ChMot}}

\newcommand{\DMgm}{\mathbf{DM}_{\mathrm{gm}}}

\newcommand{\SmProjTri}{\SmProj[\sW_{\mathrm{tri}}^{-1}]}
\newcommand{\SmProjDG}{\SmProj[\sW_{\mathrm{dg}}^{-1}]}

\newcommand{\NCSp}{\mathbf{NCSp}}
\newcommand{\NCMot}{\mathbf{NCMot}}

\renewcommand{\CH}{\mathrm{CH}}
\renewcommand{\Ho}{\mathrm{Ho}}

\makeatletter
\def\thickhrulefill{\leavevmode \leaders \hrule height 1pt\hfill \kern \z@}
\def\maketitle{%
  \thispagestyle{empty}%
  \vskip 1cm
  \begin{center}
    \Large \strut \@title \par
  \end{center}
  \par
  \begin{center}
    \emph{by } \normalfont\@author\par
  \end{center}
  }
\makeatother
\author{A. Khan Yusufzai}
\title{\uppercase{On derived categories and\\ noncommutative motives of varieties}}

\begin{document}

\maketitle

\begin{abstract}
In this short note we show how results of Orlov and To\"en imply that any equivalence between the derived categories of coherent sheaves on two varieties lifts to an equivalence at the level of dg-categories.
This establishes the link between the noncommutative geometry practised by the school of Bondal-Orlov, and the variant developed by Kontsevich and Tabuada.
As an application we recover Orlov's result that the derived category determines the Chow motive with rational coefficients up to Tate twists.
\end{abstract}

\section{Introduction}

\paratit{Overview.}

An important invariant of a scheme $X$ is its bounded derived category $\D(X)$ of coherent sheaves.
In one approach to noncommutative algebraic geometry, the class of spaces under consideration is enlarged by viewing any triangulated category as the derived category of some hypothetical space.
This direction has been studied for example by Bondal, Kuznetsov, Lunts, Orlov, Van den Bergh (e.g. \cite{orlov1997equivalences, bondal2001reconstruction, bondal2003generators, auroux2008mirror, kuznetsov2008derived}).

At the same time triangulated categories have been recognized as poorly behaved for decades: for one example let us note the lack of any kind of homotopy theory of triangulated categories.
On the other hand, in nature they arise almost always as \emph{homotopy categories} or truncations of $(\infty,1)$-categories.
M. Kontsevich has therefore suggested replacing triangulated categories by stable linear $(\infty,1)$-categories in noncommutative algebraic geometry.
Usually pretriangulated dg-categories are taken as a model for stable linear $(\infty,1)$-categories; we will recall some basic notions in section \S \ref{sec:enhancement}.

\cit{lunts2010uniqueness} showed that for a quasi-projective scheme $X$ over a commutative ring $k$, there exists a \emph{unique} dg-category $\Ddg(X)$ whose homotopy category $H^0(\Ddg(X))$ is equivalent to $\D(X)$ (we say that $\Ddg(X)$ is the unique dg-enhancement of $\D(X)$).
Hence given $X$ it is possible to consider either the triangulated category $\D(X)$ or the dg-category $\Ddg(X)$.
However it is not clear a priori whether these give the same notion of noncommutative space.
Though every dg-functor $\Ddg(X) \to \Ddg(Y)$ induces a triangulated functor $\D(X) \to \D(Y)$, the converse is not obvious: it is in general a very difficult problem to lift triangulated functors.
In particular there could be schemes $X$ and $Y$ which are isomorphic at the triangulated level, but not at the dg level.

The purpose of this short note is to show that, for smooth proper schemes $X$ over a field $k$, any fully faithful functor $\D(X) \hook \D(Y)$ lifts to a dg-functor.
This means that, in passing to the dg-version of noncommutative geometry, at least the notion of isomorphism doesn't change.
This fact is an immediate consequence of results of \cit{orlov1997equivalences} and \cit{toen2007homotopy}, but perhaps may not be obvious to non-experts.

Using the theory of noncommutative motives as developed by \cit{tabuada2011guided} we recover as an immediate corollary the facts that the derived category $\D(X)$ determines the higher algebraic K-theory, cyclic homology, and Hochschild homology of $X$.
We also recover a theorem of \cit{orlov2005derived} stating that the derived category determines the Chow motive up to Tate twists.

\paratit{Acknowledgements.}
This note is a summary of the author's master's thesis \cite{yusufzai2013perfect}.
He would like to thank his advisors, H\'el\`ene Esnault and Marc Levine, for their invaluable guidance.
He would also like to thank Dmitri Orlov and Bertrand To\"en for helpful discussions, and Jin Cao, Brad Drew and Lorenzo Mantovani for reading drafts of this note.

\paratit{Notation.}
Throughout the note, we will fix a field $k$ and work in the category $\SmProj$ of smooth projective schemes over $k$.
For $X, Y \in \SmProj$, we will write $X \times Y$ instead of $X \times_k Y$ for the fibred product over $k$.

\section{Functorial enhancement}
\label{sec:enhancement}

\para
We will use pretriangulated dg-categories as a model for stable linear $(\infty,1)$-categories.
Recall that a \dfn[]{dg-category} (\emph{differential graded category}) over $k$ is by definition a category enriched over the symmetric monoidal category $\C(k)$ of complexes of $k$-modules.
We will write $\DGCat$ for the category of small dg-categories over $k$.
Given a dg-category $\AA$, we will write $\H^0(\AA)$ for its \dfn[]{homotopy category}, given by taking the zeroth cohomologies of all the mapping complexes.
By a \dfn[]{dg-enhancement} of a triangulated category $\A$  we mean a dg-category $\AA$ whose homotopy category $\H^0(\AA)$ is equivalent to $\A$.
We refer the reader to the introductions of \cit{toen2011lectures} or \cit{keller2006differential}.

Recall that there is a natural notion of equivalence of dg-categories, called \dfn[]{quasi-equivalence}.
There is a model structure on $\DGCat$ where the weak equivalences are the quasi-equivalences (\cit{tabuada2005structure}), and we let $\DGqeq$ denote the associated homotopy category, that is to say, the localization of $\DGCat$ at the class of quasi-equivalences.

\para
For two schemes $X, Y \in \SmProj$, a \dfn[]{derived correspondence} between $X$ and $Y$ is an object of the derived category $\D(X \times Y)$.
Given derived correspondences $\sE^\bull \in \D(X \times Y)$ and $\sE'^\bull \in \D(Y \times Z)$, one defines their composite as the complex
  \[ \sE'^\bull \circ \sE^\bull = \R(p_{XZ})_*(\L(p_{XY})^*(\sE^\bull) \otimes^\L \L(p_{YZ})^*(\sE'^\bull)) \]
in $\D(X \times Z)$, where $p_{XY}$, $p_{XZ}$ and $p_{YZ}$ are the projections from $X \times Y \times Z$.
This defines a category $\DerCorr$ where morphisms are isomorphism classes of derived correspondences.

\para
Let $\TriCat$ denote the category of small triangulated categories and isomorphism classes of triangulated functors.
There is a canonical functor
  \begin{equation} \label{eq:mukai functor}
    \D : \DerCorr \to \TriCat
  \end{equation}
which maps a scheme $X \in \SmProj$ to $\D(X)$ and a derived correspondence $\sE^\bull \in \D(X \times Y)$ to the triangulated functor $\D(X) \to \D(Y)$ defined by
  \[ \D(\sE^\bull) := \R(p_Y)_*(\sE^\bull \otimes^\L \L(p_X)^*(-)), \]
where $p_X$ and $p_Y$ are the respective projections from $X \times Y$.
The functoriality of this construction was shown by \cit{mukai1981duality}.
$\D(\sE^\bull)$ is called the functor \dfn[]{represented} by $\sE^\bull$, or the \dfn[]{Fourier-Mukai functor} associated to $\sE^\bull$.

\para
By To\"en's representability theorem (\cite{toen2007homotopy}, Theorem 8.15), there are canonical bifunctorial isomorphisms of sets
  \[ \Iso(\D(X \times Y)) \isoto \Hom_{\DGqeq}(\Ddg(X), \Ddg(Y)). \]
Writing $\SmProjDG$ for the full subcategory of $\DGqeq$ spanned by the dg-categories $\Ddg(X)$ for $X \in \SmProj$, one gets a canonical equivalence of categories
  \begin{equation} \label{eq:toen functor}
    \DerCorr \isoto \SmProjDG
  \end{equation}
which is given on objects by $X \mapsto \Ddg(X)$.

\para
By Orlov's representability theorem (\cite{orlov1997equivalences}, Theorem 2.2), every fully faithful triangulated functor $\D(X) \to \D(Y)$ is represented by some derived correspondence $\sE^\bull \in \D(X \times Y)$ which is unique up to isomorphism.

By abuse of notation let $\SmProjTri \sub \TriCat$ denote the full subcategory spanned by triangulated categories of the form $\D(X)$ for $X \in \SmProj$.
Let $\SmProjTri_0 \sub \SmProjTri$ denote the (nonfull) subcategory where the morphisms are only the fully faithful functors.
We have a canonical functor
  \begin{equation} \label{eq:orlov functor}
    \SmProjTri_0 \too \DerCorr
  \end{equation}
which is \emph{faithful} (but not full).

\para
Define the \dfn[]{enhancement functor} $\varepsilon : \SmProjTri_0 \to \SmProjDG$ as the composite of the functor $\SmProjTri_0 \to \DerCorr$ (\ref{eq:orlov functor}) with the equivalence $\DerCorr \isoto \SmProjDG$ (\ref{eq:toen functor}).
  \[ \begin{tikzcd}
    \ 
      & \DerCorr \arrow{rd}{\sim}
      & \ \\
    \SmProjTri_0 \arrow{ru}\arrow{rr}{\varepsilon}
      & \ 
      & \SmProjDG
  \end{tikzcd} \]
This associates to the triangulated category $\D(X)$ its dg-enhancement $\Ddg(X)$.
Though it is not fully faithful, note that it is conservative (i.e. reflects isomorphisms).

\para
Define the \dfn[]{dehancement functor} $\delta : \SmProjDG \to \SmProjTri$ as the composite of the equivalence $\SmProjDG \isofrom \DerCorr$ with the functor $\DerCorr \to \SmProjTri$ induced by $\D : \DerCorr \to \TriCat$ (\ref{eq:mukai functor}).
  \[ \begin{tikzcd}
    \ 
      & \DerCorr \arrow{ld}[swap]{\sim}\arrow{rd}{\D}
      & \ \\
    \SmProjDG \arrow{rr}{\delta}
      & \ 
      & \SmProjTri
  \end{tikzcd} \]
This associates to the dg-category $\Ddg(X)$ its homotopy category $\D(X)$.

\para[prop:iso in PfCorr]
\begin{parathm}[Orlov, To\"en]
Let $X,Y \in \SmProj$ be smooth projective schemes over a field $k$.
The three conditions
  \begin{enumerate}
    \item $X$ and $Y$ are isomorphic in $\DerCorr$.
    \item The triangulated categories $\D(X)$ and $\D(Y)$ are equivalent.
    \item The dg-categories $\Ddg(X)$ and $\Ddg(Y)$ are quasi-equivalent.
  \end{enumerate}
are equivalent.
\end{parathm}

This is an immediate consequence of the existence of the above functors.

\section{Noncommutative motives}

\para
Let $\NCSp$ be the category of noncommutative spaces, i.e. the full subcategory of $\DGqeq$ spanned by smooth proper pretriangulated dg-categories.
The category $\NCMot$ of \dfn[]{noncommutative motives} is the karoubian envelope of the category with the same objects as $\NCSp$ and where morphisms are Grothendieck groups of internal homs; see \cit{tabuada2011guided}.
The canonical functor
  \begin{equation} \label{eq:U}
    U : \NCSp \too \NCMot
  \end{equation}
is the \dfn[]{universal additive invariant}, i.e. the universal functor sending semi-orthogonal decompositions (of the homotopy category) to direct sums (in some additive category).
See (\emph{loc. cit.}, Theorem 4.2).

\para
For $X \in \SmProj$, the dg-category $\Ddg(X)$ is smooth and proper, and hence is a noncommutative space.
In particular $\SmProjDG$ is a full subcategory of $\NCSp$.
The noncommutative motive of $X$, which we will denote $NM(X)$, is defined as the noncommutative motive of its associated noncommutative space $\Ddg(X)$.

\para[prop:NM]
\begin{paracor}
Let $X,Y \in \SmProj$ be smooth projective schemes over a field $k$.
If the triangulated categories $\D(X)$ and $\D(Y)$ are equivalent, then the noncommutative motives $NM(X)$ and $NM(Y)$ are isomorphic.
\end{paracor}

\begin{paraproof}
The enhancement functor $\varepsilon : \SmProjTri_0 \to \SmProjDG$ lifts a triangulated equivalence $\D(X) \isoto \D(Y)$ to an isomorphism $\Ddg(X) \isoto \Ddg(Y)$ in $\SmProjDG$, and therefore in $\NCSp$.
Hence $U : \NCSp \to \NCMot$ gives an isomorphism $NM(X) \isoto NM(Y)$.
\end{paraproof}

\para
Let $\Ho(\Spt)$ denote the homotopy category of spectra.
By work of \cit{blumberg2012localization}, \cit{keller1999cyclic}, \cit{schlichting2006negative} and \cit{thomason2007higher} (cf. \cit{tabuada2011guided}), each of the following can be defined as functors $\NCSp \to \Ho(\Spt)$ and are additive invariants in the above sense:
  \begin{enumerate}
    \item algebraic K-theory,
    \item cyclic homology,
    \item topological cyclic homology,
    \item Hochschild homology,
    \item topological Hochschild homology.
  \end{enumerate}
By universality of the functor $U : \NCSp \to \NCMot$ (\ref{eq:U}), one has the following corollary.

\begin{paracor}
Let $X,Y \in \SmProj$ be smooth projective schemes over a field $k$ and let $H_*$ denote one of the above functors.
If the triangulated categories $\D(X) \shiso \D(Y)$ are equivalent, then $H_*(X)$ and $H_*(Y)$ are isomorphic.
\end{paracor}

\section{Chow motives}

\para[lem:orbit]
Let $\A$ be an additive category and let $T : \A \isoto \A$ be an auto-equivalence.
The \dfn[]{orbit category of $\A$ with respect to $T$} is the category $\A/T$ whose objects are the same as those of $\A$ and whose morphisms are given by
  \[ \Hom_{\A/T}(X, Y) = \bigoplus_{i \in \Z} \Hom_\A(X, T^i(Y)) \]
for all $X,Y \in \A$.
The law of composition is defined as follows: for two morphisms $f = (f^i)_i : X \to Y$ and $g = (g^j)_j : Y \to Z$ in $\A/T$, the composite $g \circ f$ is defined as the morphism whose $k$-th component is the sum
  \[ (g \circ f)^k = \sum_{i+j=k} T^i(g^j) \circ f^i. \]
Let $\pi_{\A/T} : \A \to \A/T$ denote the canonical functor which maps a morphism $f$ to the morphism which is $f$ in the zeroth component and $0$ everywhere else.

\begin{paralemma}
Let $\A$ be an additive category and $T : \A \isoto \A$ an additive auto-equivalence.
Suppose that $\A$ admits arbitrary direct sums.
Then the projection functor $\pi = \pi_{\A/T} : \A \to \A/T$ admits a right adjoint
  \[ \tau = \tau_{\A/T} : \A/T \too \A \]
which maps an object $X$ to the direct sum of all the objects $T^i(X)$ $(i \in \Z$).
\end{paralemma}

\para
Let $\ChMot(\Q)$ be the category of Chow motives with rational coefficients, and $M : \SmProj^\op \to \ChMot(\Q)$ the canonical functor (see \cit{andre2004introduction}).
Let $\Q(1) \in \ChMot(\Q)$ denote the Tate motive, and recall that the functor $\Q(1) \otimes - : \ChMot(\Q) \to \ChMot(\Q)$ is an auto-equivalence.
Let $\ChMot(\Q)/\Q(1)$ denote the associated orbit category; we call this the \dfn[]{category of Chow motives modulo Tate twists}.
We will write $M(i) = M \otimes \Q(i) = M \otimes \Q(1)^{\otimes i}$ for a motive $M \in \ChMot(\Q)$.

Let $\DMgm(\Q)$ be the triangulated category of geometric motives of Voevodsky and recall that there is a canonical functor
  \begin{equation}
    \ChMot(\Q) \too \DMgm(\Q) \label{eq:chm to dm}
  \end{equation}
which sends the Tate motive $\Q(1)$ to $\Q(-1)[-2]$ (where $-[n]$ denotes the $n$-fold composition of the translation functor); see (\emph{loc. cit.}, Th\'eor\`eme 18.3.1.1).
Hence one gets an induced functor $\ChMot(\Q)/\Q(1) \too \DMgm(\Q)/\Q(-1)[-2]$ on the orbit categories and by the above lemma (\ref{lem:orbit}), as the triangulated category $\DMgm(\Q)$ admits arbitrary direct sums, one gets a canonical functor
  \begin{equation}
    \ChMot(\Q)/\Q(1) \too \DMgm(\Q)/\Q(-1)[-2] \stackrel{\tau}{\too} \DMgm(\Q) \label{eq:orbit}
  \end{equation}
which maps a motive to the direct sum of all its Tate twists:
  \[ M \mapsto \bigoplus_{i \in \Z} M(i)[2i]. \]

\para
\cit{kontsevich2009notes} noted that there is a canonical fully faithful functor
  \[ \nu : \ChMot(\Q)/\Q(1) \hooklong \NCMot(\Q) \]
which is given on morphisms by
  \[ \bigoplus_{i} \CH^i(X \times Y, \Q) \isoto K_0(X \times Y) \otimes \Q, \]
the inverses of the Grothendieck-Riemann-Roch isomorphisms $\ch(-) \cdot \sqrt{\td_{X \times Y}}$.
Note that this fits into a commutative diagram
  \[ \begin{tikzcd}
    \SmProjDG \arrow[hook]{r}\arrow{d}[swap]{M}
      & \NCSp \arrow{d}{U} \\
    \ChMot(\Q)/\Q(1) \arrow[hook]{r}{\nu}
      & \NCMot(\Q).
  \end{tikzcd} \]

\para
\begin{parathm}[\cit{orlov2005derived}]
Let $X, Y \in \SmProj$ be smooth projective schemes over a field $k$.
If the triangulated categories $\D(X) \shiso \D(Y)$ are equivalent, then there is an isomorphism
  \[ \bigoplus_{i\in\Z} M(X)(i)[2i] \iso \bigoplus_{j\in\Z} M(Y)(j)[2j] \]
in the triangulated category $\DMgm(\Q)$ of geometric motives.
\end{parathm}

\begin{paraproof}
Since $NM(X) \shiso NM(Y)$ in $\NCMot(\Q)$ by (\ref{prop:NM}), and the functor $\nu : \ChMot(\Q)/\Q(1) \hook \NCMot(\Q)$ is fully faithful, one has $M(X) \shiso M(Y)$ in $\ChMot(\Q)/\Q(1)$.
Then the functor $\ChMot(\Q)/\Q(1) \to \DMgm(\Q)$ (\ref{eq:orbit}) gives the desired isomorphism.
\end{paraproof}

\bibliographystyle{plainnat-aky}

{\small
\setlength{\bibsep}{0.0pt}
\bibliography{references}

\begin{thebibliography}{20}
\providecommand{\natexlab}[1]{#1}
\providecommand{\url}[1]{\texttt{#1}}
\expandafter\ifx\csname urlstyle\endcsname\relax
  \providecommand{\doi}[1]{doi: #1}\else
  \providecommand{\doi}{doi: \begingroup \urlstyle{rm}\Url}\fi

\bibitem[Andr{\'e}(2004)]{andre2004introduction}
Yves Andr{\'e}.
\newblock \emph{Une introduction aux motifs: motifs purs, motifs mixtes,
  p{\'e}riodes}.
\newblock Panoramas et Synth{\`e}ses. Soci{\'e}t{\'e} Math{\'e}matique de
  France, 2004.
\newblock ISBN 9782856291641.

\bibitem[Auroux et~al.(2008)Auroux-Katzarkov--Orlov]{auroux2008mirror}
Denis Auroux, Ludmil Katzarkov, and Dmitri Orlov.
\newblock Mirror symmetry for weighted projective planes and their
  noncommutative deformations.
\newblock \emph{Ann. of Math.(2)}, 167\penalty0 (3):\penalty0 867--943, 2008.

\bibitem[Blumberg-Mandell(2012)]{blumberg2012localization}
Andrew~J. Blumberg and Michael~A. Mandell.
\newblock Localization theorems in topological {H}ochschild homology and
  topological cyclic homology.
\newblock \emph{Geometry \& Topology}, 16\penalty0 (2):\penalty0 1053--1120,
  2012.

\bibitem[Bondal-Orlov(2001)]{bondal2001reconstruction}
Alexei Bondal and Dmitri Orlov.
\newblock Reconstruction of a variety from the derived category and groups of
  autoequivalences.
\newblock \emph{Compositio Mathematica}, 125\penalty0 (03):\penalty0 327--344,
  2001.

\bibitem[Bondal-Van~den Bergh(2003)]{bondal2003generators}
Alexei Bondal and Michel Van~den Bergh.
\newblock Generators and representability of functors in commutative and
  noncommutative geometry.
\newblock \emph{Mosc. Math. J}, 3\penalty0 (1):\penalty0 1--36, 2003.

\bibitem[Keller(1999)]{keller1999cyclic}
Bernhard Keller.
\newblock On the cyclic homology of exact categories.
\newblock \emph{Journal of Pure and Applied Algebra}, 136\penalty0
  (1):\penalty0 1--56, 1999.

\bibitem[Keller(2006)]{keller2006differential}
Bernhard Keller.
\newblock On differential graded categories.
\newblock \emph{arXiv preprint math/0601185}, 2006.

\bibitem[Kontsevich(2009)]{kontsevich2009notes}
Maxim Kontsevich.
\newblock Notes on motives in finite characteristic.
\newblock In \emph{Algebra, Arithmetic, and Geometry}, pages 213--247.
  Springer, 2009.

\bibitem[Kuznetsov(2008)]{kuznetsov2008derived}
Alexander Kuznetsov.
\newblock Derived categories of quadric fibrations and intersections of
  quadrics.
\newblock \emph{Advances in Mathematics}, 218\penalty0 (5):\penalty0
  1340--1369, 2008.

\bibitem[Lunts-Orlov(2010)]{lunts2010uniqueness}
Valery Lunts and Dmitri Orlov.
\newblock Uniqueness of enhancement for triangulated categories.
\newblock \emph{Journal of the American Mathematical Society}, 23\penalty0
  (3):\penalty0 853--908, 2010.

\bibitem[Mukai(1981)]{mukai1981duality}
Shigeru Mukai.
\newblock Duality between {$D(X)$} and {$D(\hat{X})$} with its application to
  {P}icard sheaves.
\newblock \emph{Nagoya Math. J.}, 81:\penalty0 153--175, 1981.

\bibitem[Orlov(1997)]{orlov1997equivalences}
Dmitri Orlov.
\newblock Equivalences of derived categories and {K}3 surfaces.
\newblock \emph{Journal of Mathematical Sciences}, 84\penalty0 (5):\penalty0
  1361--1381, 1997.

\bibitem[Orlov(2005)]{orlov2005derived}
Dmitri Orlov.
\newblock Derived categories of coherent sheaves and motives.
\newblock \emph{Russian Mathematical Surveys}, 60\penalty0 (6):\penalty0
  1242--1244, 2005.

\bibitem[Schlichting(2006)]{schlichting2006negative}
Marco Schlichting.
\newblock Negative {K}-theory of derived categories.
\newblock \emph{Mathematische Zeitschrift}, 253\penalty0 (1):\penalty0 97--134,
  2006.

\bibitem[Tabuada(2005)]{tabuada2005structure}
Gon\c{c}alo Tabuada.
\newblock Une structure de cat{\'e}gorie de modeles de {Q}uillen sur la
  cat{\'e}gorie des dg-cat{\'e}gories.
\newblock \emph{Comptes Rendus Mathematique}, 340\penalty0 (1):\penalty0
  15--19, 2005.

\bibitem[Tabuada(2011)]{tabuada2011guided}
Gon\c{c}alo Tabuada.
\newblock A guided tour through the garden of noncommutative motives.
\newblock \emph{arXiv preprint arXiv:1108.3787}, 2011.

\bibitem[Thomason-Trobaugh(2007)]{thomason2007higher}
Robert~W. Thomason and Thomas Trobaugh.
\newblock Higher algebraic {K}-theory of schemes and of derived categories.
\newblock In \emph{The Grothendieck Festschrift Volume III}, pages 247--435.
  Springer, 2007.

\bibitem[To{\"e}n(2007)]{toen2007homotopy}
Bertrand To{\"e}n.
\newblock The homotopy theory of dg-categories and derived {M}orita theory.
\newblock \emph{Inventiones mathematicae}, 167\penalty0 (3):\penalty0 615--667,
  2007.

\bibitem[To{\"e}n(2011)]{toen2011lectures}
Bertrand To{\"e}n.
\newblock Lectures on dg-categories.
\newblock \emph{Topics in algebraic and topological K-theory}, pages 243--301,
  2011.

\bibitem[Yusufzai(2013)]{yusufzai2013perfect}
A.~Khan Yusufzai.
\newblock Perfect correspondences and {C}how motives.
\newblock \emph{arXiv preprint arXiv:1310.0249}, 2013.

\end{thebibliography}
}

\vspace{1.2em}

{\small
\noindent
Fakultät Mathematik,
Universität Duisburg-Essen,
Thea-Leymann-Straße 9,
45127 Essen,
Germany\\
\emph{E-mail address:} \texttt{khan.adeel@stud.uni-due.de}\\
\emph{URL:} \texttt{http://www.preschema.com}}

\end{document}